# A NOTE ON THE ENUMERATION OF DIRECTED ANIMALS VIA GAS CONSIDERATIONS

By Marie Albenque

*Université Paris Diderot—Paris 7*

In the literature, most of the results about the enumeration of directed animals on lattices via gas considerations are obtained by a formal passage to the limit of enumeration of directed animals on cyclical versions of the lattice.

Here we provide a new point of view on this phenomenon. Using the gas construction given in [*Electron. J. Combin.* (2007) **14** R71], we describe the gas process on the cyclical versions of the lattices as a cyclical Markov chain (roughly speaking, Markov chains conditioned to come back to their starting point). Then we introduce a notion of convergence of graphs, such that if $(G_n) \to G$ then the gas process built on $G_n$ converges in distribution to the gas process on $G$. That gives a general tool to show that gas processes related to animals enumeration are often Markovian on lines extracted from lattices.

We provide examples and computations of new generating functions for directed animals with various sources on the triangular lattice, on the $\mathcal{T}_n$ lattices introduced in [*Ann. Comb.* **4** (2000) 269–284] and on a generalization of the $\mathcal{L}_n$ lattices introduced in [*J. Phys. A* **29** (1996) 3357–3365].

**1. Introduction.** Let $G = (V, E)$ be a *directed* graph with set of vertices $V$ and set of oriented edges $E$. Let $A$ and $S$ be two subsets of $V$, with $S \subset A$. We say that $A$ is a *directed animal* (DA) with source $S$ if and only if any vertex of $A$ can be reached from an element of $S$ through a directed path having all its vertices in $A$ (see Figure 1). The vertices of $A$ are called *cells* and the number of cells, denoted $|A|$, is the *area* of $A$. We denote $\mathcal{G}_S^G$ the generating function (GF) for DA on $G$ with source $S$ counted according to

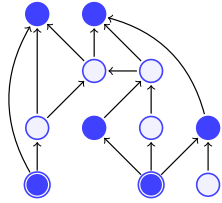

Fig. 1. *Example of a DA with area 6. The cells of the DA are dark and the vertices of the source double circled.*

their area

$$\mathcal{G}_S^G(t) = \sum_{\substack{A, \text{ DA} \\ \text{with source S}}} t^{|A|} = \sum_{k \geq |S|} a_k t^k,$$

where $a_k$ is the number of DA on $G$ with source $S$ and area $k$.

In the following, we will always assume that the cells of $S$ form an independent set on the directed graph $G$—we say that $S$ is a free set—the formal definition follows.

DEFINITION 1.  Let $G = (V, E)$ be an oriented graph and $x$ and $y$ be two vertices of $G$. We say that $x$ is a father of $y$ or equivalently that $y$ is a child of $x$ if there is an edge from $x$ to $y$.

More generally, $x$ is called an *ancestor* of $y$ if there exists a directed path from $x$ to $y$.

Let now $S$ be a subset of $V$; we say that $S$ is a *free set* of vertices of $G$ if and only if for every $x, y \in S$ such that $x \neq y$, $x$ is not an ancestor of $y$.

In this article we focus on the link between enumeration of DA and hard particle gas models.

DEFINITION 2.  Let $G = (V, E)$ be a graph, a *gas occupation* or *gas configuration* on $G$ is a map $X$ from $V$ to $\{0, 1\}$. The vertices $v \in V$ such that $X(v) = 1$ are said to be occupied, the others are said to be empty.

A *hard particle* gas occupation of a graph is a gas occupation with the additional constraint that two occupied vertices cannot be neighbors (the occupied cells form then an independent set).

A *gas model* is a probability law on gas occupations. For a given gas model, we call *density* in a vertex $v$ the probability for $v$ to be occupied, that is $\mathbb{P}(X(v) = 1)$.

Since the pioneering work of Dhar [7], the connection between DA and gas models have been widely exploited. We shall now give a short overview of the different contributions on this subject (we refer the reader to [4] and [10] for more exhaustive references). In Dhar [7], using some statistical mechanics shows that computing the area generating function for DA on the square lattice is equivalent to computing the density of a hard particle gas model. This result was obtained after Nadal, Derrida and Vannimenus [11] and Hakim and Nadal [9] obtained the generating function of DA on some "cylindric" square lattices.

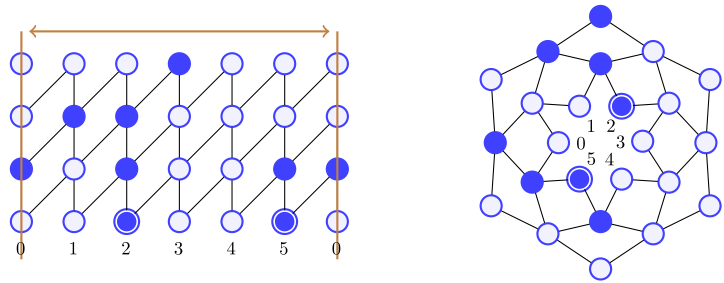

FIG. 2. *The same DA on two representations of the cylindric version of the square lattice with a width of 6.*

Those "cylindric" lattices are defined as follows: Let $G$ be an oriented lattice—that is, an oriented translation-invariant graph—with its vertices indexed by a subset of $\mathbb{Z}^2$. If we consider that the abscissa of the vertices of $G$ are labeled by elements of $[N] := \mathbb{Z}/N\mathbb{Z}$ instead of $\mathbb{Z}$, we obtain the width-bounded variant of $G$ with cyclic boundary conditions (see Figure 2). We denote it $G^{(N)}$ and call it the cyclic or cylindric version of $G$ of width $N$.

Bousquet-Mélou [4] extends Dhar's correspondence between the hard particle gas models and enumeration of DA on cyclic square lattices. Particularly, she shows that gas models allow the enumeration of DA not only according to their area but also for instance according to their left perimeter or their number of loops. Those results were then generalized to a family of lattices in a joint work with Conway [5]. In Conway [5] and Bousquet-Mélou [4], the gas models studied are defined on the cylindric versions of graphs and the GF for DA is obtained as the formal limit of the density of the gas when the width grows to infinity. Since computing the density of the gas model is not always tractable, the former result does not necessarily lead to effective results about enumeration of DA. However, that new link establishes gas models as a powerful and polyvalent tool for the counting of DA.

In the latter works, the link between DA and gas is formal and appears because DA and gas models are shown to verify the same recursive decomposition along with the layers of the graph. It notably implies that that approach is only valid for graphs that can be decomposed nicely into layers.

Le Borgne and Marckert [10] give a new insight into the connection between gas and DA. They construct a coupling between random DA and random gas models and give a combinatorial proof that for a free set $S$ the GF of DA with source $S$ is equal to the probability for the vertices of $S$ to be occupied (a construction of that coupling is sketched in Section 2). Contrary to the construction on cylinders, in [10] the gas model is well defined on any acyclic graph and in particular on the whole lattice, where some of its stochastic properties can be studied. On the square lattice for instance, its restriction to a line is shown to be Markovian, which allows to compute explicitly the GF for DA with any source included in a line.

We must mention that there exist other fruitful approaches to the combinatorics of DA. Some results have been obtained by establishing links with heaps of pieces introduced by Viennot in [12] (see for instance [2, 3, 6, 13]) or with paths in the plane [8] or via the ECO method [1].

We now describe the content of this paper and its organization. Our aim here is to give a general framework that allows to reduce the enumeration of DA with various sources on a graph $G$ to the same enumeration on "simpler" versions of $G$. As mentioned above, simplifying the graphs we work on is a classical idea. Here, the difference with the works cited above relies on the fact that thanks to the gas construction given in [10], we can now study the convergence of the gas models as stochastic processes and not only the formal convergence of their density. This leads both to a better understanding of the gas models and to new results about enumeration of DA with various sources.

The first point is to make the notion of "simpler" versions of $G$ accurate; in Section 3.1 we provide a distance on the set of graphs with marked vertices, corresponding to sources [see equation (3.1)]. Roughly speaking, $G_n$ converges to $G$ for the notion of convergence of graphs induced by that distance (which corresponds roughly to the convergence of the neighborhood of sources) implies that $\mathcal{G}_S^{G_n}$ converges to $\mathcal{G}_S^G$. In terms of probability, that means the convergence of the finite-dimensional laws of the gas under additional assumptions (Theorem 2).

Then we need to compute the law of the limiting gas process obtained thanks to that convergence. That is possible on some lattices. The multiplicative formula obtained for the distribution of gas restricted to a line on the cylinder in [4] and [5] leads to the intuition that that multiplicative structure may be preserved when the width of the cylinder goes to infinity and that the limiting process obtained above should be Markovian. For that reason, in Section 3.2 we define a *cyclic Markov chain* as a Markov chain conditioned to come back to its initial state after a fixed number of steps (Definition 6). We then give a representation of the gas on the cylinder as a cyclic Markov chain. Then in broad terms when the width of the cylinder grows, the conditioning induced by the cyclic condition is less and less constraining. At the limit, it eventually disappears which therefore yields that the limiting process is Markovian. We provide in Theorem 3 a formal statement of those two ideas; that provides a frame in which the gas process on a line is Markovian.

We define in Section 4.1 the family of lattices $(\mathcal{L}_R)_{R \subset \mathbb{N}}$, which extends the family of lattices $(\mathcal{L}_n)_{n \geq 2}$ introduced in [5]. We apply Theorem 2 and Theorem 3 to it, to the triangular lattice and to the family of $\mathcal{T}_n$ lattices introduced in [6]. In Section 4, we show that for those three examples, the restriction of the gas process to a line is Markovian. Thanks to the link between gas models and GF of DA, that allows us to obtain some GF for DA with various sources; see, for example, Proposition 1 for some results on the triangular lattice.

**2. Definition of the gas model.** We sketch the construction of the gas model given in [10] and its link with enumeration of DA according to their area. Let $G = (V, E)$ be a directed graph without multiple edges nor directed cycles and such that the number of children of each node is finite.

The probability space we work on is $\Omega = \{a, b\}^V$ endowed with the $\sigma$-field generated by the finite subsets of vertices. We equip that space with the product probability $\mathbb{P}_p = (p\delta_a + (1-p)\delta_b)^{\otimes V}$, where $\delta_a$ is the standard Dirac measure on $\{a\}$. In other terms, $\omega \in \Omega$ is a coloring of $G$ and under $\mathbb{P}_p$ each vertex has, independently of the others, color $a$ or $b$ with respective probabilities $p$ and $1-p$. For $x \in V$, $\omega(x)$ gives the color of $x$. From that random coloring we construct DA and a model of gas. Notice that the DA and gas process defined below are deterministic functions of the random coloring.

DEFINITION 3. Let $S$ be a subset of $V$ and $\omega$ be a random coloring of $G$. We denote by $S_\bullet(\omega) = \{x \in S, \omega(x) = a\}$, the (random) subset of $S$ with color $a$. We then define the random variable $\mathbf{A}^S$ as the maximal DA for the inclusion partial order with source $S_\bullet(\omega)$ and set of cells the $a$-colored vertices $x$ that can be reached from $S_\bullet(\omega)$ by an $a$-colored path (see Figure 3).

For a set $S$ such that $|S| \geq 1$, the random DA $\mathbf{A}^S$ may be infinite with positive probability. Let $p_{\mathrm{crit}}^G$ be the threshold for the existence of an infinite DA with positive probability (it corresponds to the critical probability for the oriented percolation on $G$)

$$(2.1) \qquad p_{\mathrm{crit}}^G := \sup\{p, S : \mathbb{P}_p(|\mathbf{A}^S| < \infty) = 1 \text{ and } |S| < \infty\}.$$

For a general graph $G$, $p_{\mathrm{crit}}^G$ is difficult to compute and can even be equal to zero. In the examples given in Section 4, the outdegree of any vertex is bounded and in that case $p_{\mathrm{crit}} > 0$ (see, e.g., Proposition 2.2 of [10]). For any $p < p_{\mathrm{crit}}^G$ a gas occupation $X^G$ on $G$ is defined from a random coloring $\omega = (\omega(v))_{v \in V}$ as follows (see Figure 3 for an example):

$$(2.2) \qquad X^G(v) = \begin{cases} 0, & \text{if } \omega(v) = b, \\ \displaystyle\prod_{v' \text{children of } v} (1 - X^G(v')), & \text{if } \omega(v) = a. \end{cases}$$

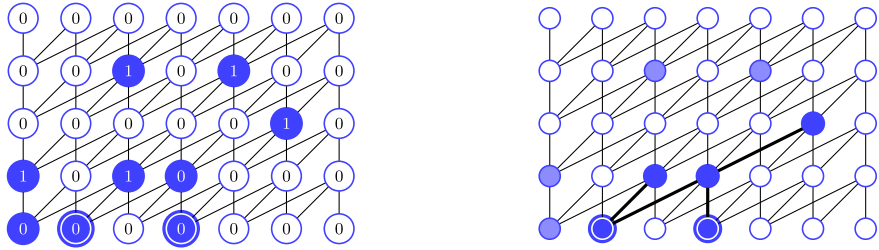

Fig. 3. *The gas occupation (on the left) and the DA $\mathbf{A}^S$ on $\mathcal{L}_3$, obtained from the same coloring of the vertices. Cells colored with a (respectively b) are dark (respectively white) and the vertices of S are double circled.*

The definition of $p_{\mathrm{crit}}^G$ ensures that the gas process is almost surely well defined as its recursive computation ends within a finite number of steps for any $p < p_{\mathrm{crit}}^G$ (see Proposition 2.4 of [10] for details). From now on we always assume that $p < p_{\mathrm{crit}}^G$ and that the gas model considered is the probability law denoted $\mathbb{P}_p^G$ induced by that construction.

The link between enumeration of DA and that gas model is given by the following result:

THEOREM 1 (Le Borgne and Marckert [10]). *Let $G = (V, E)$ be a directed graph and $S$ be a free set of $G$. For any $p$ in $[0, R_S^G)$, we have*

$$(2.3) \qquad \mathbb{P}_p^G(X^G(v) = 1, v \in S) = (-1)^{|S|} \mathcal{G}_S^G(-p),$$

*where $R_S^G$ is the radius of convergence of $\mathcal{G}_S^G$.*

With that theorem, the computation of the generating function for DA comes down to the computation of the probability for some vertices to be occupied for the gas model $\mathbb{P}_p^G$. That explains why in the next section we focus only on the study of the gas model and resume the enumeration of DA in Section 4.

**3. Convergence of graphs, gas models and DA.** We develop in that section some tools allowing us to reduce to simpler graphs the study of the stochastic properties of a gas model on a graph.

3.1. *Convergence of graphs.* As recalled in the Introduction, most of the results obtained about the enumeration of DA on a lattice $G$ via the study of gas models have been proved by a passage to the limit. More precisely, the gas models are studied on $G^{(N)}$, the cylindric version of $G$ (see the Introduction). For a fixed size $n$, the set of DA with size $n$ coincide on $G$ and $G^{(N)}$ when say $N \gg n$. It amounts to saying that

$$\mathcal{G}_{\{x\}}^{G^{(N)}} \underset{N}{\longrightarrow} \mathcal{G}_{\{x\}}^{G}$$

or equivalently in the gas model's point of view that the density of the gas converges formally (in the sense that $\sum a_{n,k} x^k \longrightarrow_n \sum a_k x^k$ if and only if $a_{n,k} \longrightarrow_n a_k$ for every $n \in \mathbb{N}$).

The aim of this section is to make clear a notion of convergence of graphs (i.e., a topology on the set of graphs) which induces the convergence of the finite-dimensional laws of the gas process and hence the convergence of the generating function of DA. That convergence is no longer seen only as a formal convergence of generating functions but as the convergence of the distribution of a stochastic process.

In the following, we always assume that the graphs considered are directed, without directed cycles nor multiple edges and that the number of children of each node is finite (a node can, though, have an infinite number of parents) so that the gas model given in Section 2 is defined.

DEFINITION 4. We call *marked directed graph*, a pair $(G = (V,E), Z)$ where $Z$ is a subset of $V$. We denote by $V_Z$ the subset of $V$ of nodes having at least one ancestor in $Z$, and by $G(Z)$ the subgraph of $G$ having as set of nodes $V_Z$ (and set of edges the edges of $E$ linking them).

To see $Z$ as a source and $G(Z)$ as the maximal DA on $G$ with source $Z$ may help to understand better Theorem 2.

DEFINITION 5. Two directed marked graphs $(G = (V,E), Z)$ and $(G' = (V',E'), Z')$ are said to be *isomorphic*—we write $(G, Z) \sim (G', Z')$—if $G(Z)$ and $G'(Z')$ are equal up to a relabeling of the vertices, in other words, if there exists a bijective application $\phi$ from $V_Z$ onto $V'_{Z'}$ such that for any $x, y$ in $V_Z$, $(x, y) \in E$ is equivalent to $(\phi(x), \phi(y)) \in E'$.

The relation $\sim$ is an equivalence relation on the set of marked directed graphs. We denote by $\mathcal{O}$ the set of directed graph quotiented by that relation. For any marked graph $(G, Z)$ we denote by $\overline{(G, Z)}$ its class in $\mathcal{O}$.

We denote $\mathcal{A}_Z^G$ the set of DA on $G$ with source $Z$. The graph $(G, Z)$ is the right (or minimal) structure that provides all the knowledge necessary to study the gas configuration on $Z$ and the DA with source $Z$ in $G$ [that depends also on the coloring on $(G, Z)$]. From the construction of the gas
model and random DA given in Section 2, it is clear that if $(G, Z) \sim (G', Z')$
then $|Z| = |Z'|$ and $\mathcal{G}_Z^G = \mathcal{G}_{Z'}^{G'}$ and the application $\phi$ provides a probability
isomorphism between the gas occupations on $Z$ and $Z'$, which implies that
$\mathbb{P}_p^G(X_s^G = 1, s \in Z) = \mathbb{P}_p^{G'}(X_s^{G'} = 1, s \in Z')$.

For any $r \geq 0$, we define $B_r(G, Z)$ as the subgraph of $(G, Z)$ containing
only the vertices $v$ of $(G, Z)$ such that $d(v, Z) = \inf_{u \in Z} d(u, v) \leq r$, where
the distance must be understood as a directed distance on graphs, that is,

$$d(u, v) = \inf\{|w|, \text{ where } w \text{ is an } \textit{oriented} \text{ path from } u \text{ to } v\}.$$

As announced above, we now define a distance $d_{\mathcal{O}}$ on $\mathcal{O}$ which gives a
suitable notion of convergence of graphs: For any $O$ and $O'$ in $\mathcal{O}$, we set

$$(3.1) \qquad d_{\mathcal{O}}(O, O') = \inf\left\{\frac{1}{r+1}, r \text{ such that } B_r(G, Z) \sim B_r(G', Z')\right\},$$

where $(G, Z) \in O$ and $(G', Z') \in O'$ [we let the reader check that that is
indeed a distance in $\mathcal{O}$ and, in particular, that it does not depend on the
choices of $(G, Z)$ and $(G', Z')$].

THEOREM 2.   *Let* $(G_n = (V_n, E_n), Z_n)$ *be a sequence of directed marked
graphs, and* $(G = (V, E), Z)$ *be a directed marked graph. Let* $a_{n,k} = \#\{A \in \mathcal{A}_{Z_n}^{G_n}, |A| = k\}$ *be the number of DA with source* $Z_n$ *in* $G_n$ *having* $k$ *cells, and
denote by* $a_k = \#\{A \in \mathcal{A}_Z^G, |A| = k\}$.
*If* $d_{\mathcal{O}}((\overline{G_n, Z_n}), (\overline{G, Z})) \to 0$ *then:*

1. $\mathcal{G}_{Z_n}^{G_n}(p) = \sum_{k \geq |Z_n|} a_{n,k} p^k \xrightarrow{}_{n \to \infty} \mathcal{G}_Z^G(p) = \sum_{k \geq |Z|} a_k p^k$ *where the conver-
   gence holds formally in the set of formal series with coefficient in* $\mathbb{N}$ *(i.e.,
   for any* $k$, $a_{n,k} \to a_k$ *when* $n \to \infty$).
2. *If there exists* $c, d \geq 0$ *such that for any* $n$ *large enough,*

$$(3.2) \qquad\qquad a_{n,k} \leq cd^k \qquad \text{for any } k \geq 1,$$

   *then for any* $p < 1/d$, *the finite-dimensional laws of the gas occupation on
   $Z_n$ according to* $\mathbb{P}_p^{G_n}$ *converge towards those on* $Z$ *distributed according
   to* $\mathbb{P}_p^G$, *that is,*

$$\mathbb{P}_p^{G_n}(X_s^{G_n} = 1, s \in Z_n) \xrightarrow{}_{n} \mathbb{P}_p^G(X_s^G = 1, s \in Z).$$

PROOF.   1. First, if $d_{\mathcal{O}}((\overline{G_n, Z_n}), (\overline{G, Z})) \to 0$, then for any $r$, when $n$ is
large enough, the two graphs $B_r(G_n, Z_n)$ and $B_r(G, Z)$ are isomorphic. That
implies that the coefficients of $\mathcal{G}_{Z_n}^{G_n}$ and $\mathcal{G}_Z^G$ coincide at least up to the $r$th.

2. First, condition (3.2) implies that $p_{\text{crit}}^{G_n} \geq 1/d$, therefore the gas model
$\mathbb{P}_p^{G_n}$ is well defined for any $p < 1/d$.

From the construction of the gas model, we can notice that the event $\{X_s^G = 1, s \in Z\}$ does not depend on the coloring of all the vertices of $G$ but only on vertices of $\mathbf{A}^Z$ (see Definition 3). Since we assume $p < 1/d$, $\mathbf{A}^Z$ is almost surely finite according to $\mathbb{P}_p^G$; that implies that for any $\varepsilon > 0$, there exists $m_\varepsilon$ such that $\mathbb{P}_p^G(|A^Z| \geq m_\varepsilon) < \varepsilon$.

As when $n$ is large enough, the two graphs $B_{m_\varepsilon}(G_n, Z_n)$ and $B_{m_\varepsilon}(G, Z)$ are isomorphic, there exists an application $\phi$ that maps $B_{m_\varepsilon}(G_n, Z_n)$ onto $B_{m_\varepsilon}(G, Z)$. Thus $\phi$ induces a probability isomorphism between the coloring of $B_{m_\varepsilon}(G_n, Z_n)$ and of $B_{m_\varepsilon}(G, Z)$. Therefore, conditionally on the event $\{|\mathbf{A}^Z| < m_\varepsilon\}$, the image of $\mathbf{A}^{Z_n}$ by $\phi$ is $\mathbf{A}^Z$ and we get

$$\mathbb{P}_p^{G_n}(X_{s}^{G_n} = 1, s \in Z_n || A^Z| < m_\varepsilon) = \mathbb{P}_p^G(X_s^G = 1, s \in Z || A^Z| < m_\varepsilon).$$

This concludes the proof, since $\mathbb{P}_p^G(|\mathbf{A}^Z| < m_\varepsilon) \geq 1 - \varepsilon$ by definition of $m_\varepsilon$. $\square$

REMARK 1. Even if in the applications of the latter theorem in Section 4, we always assume that $Z_n$ and $Z$ are free sets. There is no such assumption in the theorem and $Z_n$ and $Z$ can be any sets.

### 3.2. A variation on Markov processes.

The spirit of this section is guided by the results obtained for enumeration of DA in [4] and [5]. It often happens that the probability distribution of the gas has a multiplicative form on cylinders. That leads to the intuition that the limiting process obtained when the width goes to infinity is Markovian. We give here an appropriate frame to make that intuition rigorous.

In this section, we always assume that $E$ is a finite state space, $\nu$ a probability measure on $E$ and $\mathbf{M}$ a stochastic matrix on $E$. We say that $Y = (Y_i)_{i \in \mathbb{N}}$ is a $(\nu, \mathbf{M})$-MC if it is a Markov chain with $\nu$ as initial law and $\mathbf{M}$ as transition matrix.

DEFINITION 6. For any nonnegative $N$, we call *cyclic Markov chain of length $N$ on $E$* with initial law $\nu$ and transition matrix $\mathbf{M}$, a process $(X_i)_{i \in \{0,...,N-1\}}$ which is a Markov chain conditioned to come back to its starting point after $N$ steps and we say that $X$ is a $(\nu, \mathbf{M}, N)$-cyclic MC.

Let $Y$ be a $(\nu, \mathbf{M})$-MC, for any $x_0, \ldots, x_{N-1} \in E$, the law of $(X_i)_{i \in \{0,...,N-1\}}$ is equal to

$$(3.3) \quad \begin{aligned} &\mathbb{P}(X_0 = x_0, \ldots, X_{N-1} = x_{N-1}) \\ &= \mathbb{P}(Y_0 = x_0, \ldots, Y_{N-1} = x_{N-1} | Y_0 = Y_N). \end{aligned}$$

In other words,

$$(3.4) \quad \mathbb{P}(X_0 = x_0, \ldots, X_{N-1} = x_{N-1}) = \frac{\nu(x_0) \prod_{i=0}^{N-1} \mathbf{M}_{x_i, x_{i+1}}}{\widetilde{Z_N}},$$

where $x_N = x_0$ and $\widetilde{Z_N} = \sum_{x'_0,\dots,x'_{N-1}} \nu(x'_0) \prod_{i=0}^{N-1} \mathbf{M}_{x'_i,x'_{i+1}}$.

Note that if $X$ is a $(\nu, \mathbf{M}, N)$-cyclic MC, the distribution of $X_0$ is given by

$$(3.5) \qquad \mathbb{P}(X_0 = x) = \frac{\nu(x)(\mathbf{M}^N)_{x,x}}{\widetilde{Z_N}} \qquad \text{for any } x \in E$$

and the distribution of $X_1$ by

$$(3.6) \qquad \mathbb{P}(X_1 = x_1) = \left( \sum_{x_0} \nu(x_0) \mathbf{M}_{x_0,x_1} (\mathbf{M}^{N-1})_{x_1,x_0} \right) (\widetilde{Z_N})^{-1}.$$

Equation (3.5) implies that the distribution of $X_0$ is not $\nu$ except for exceptional cases. Combining equations (3.5) and (3.6) implies that if $\nu = \mathcal{U}_E$, the uniform law on $E$, then the cyclic MC is stationary, that is, for any $x \in E$, $\mathbb{P}(X_i = x) = \mathbb{P}(X_0 = x)$.

On the other hand, assume that the initial law $\nu$ is an invariant law for $\mathbf{M}$, then a $(\nu, \mathbf{M}, N)$-cyclic MC is not necessarily stationary. Roughly speaking, the term $(\mathbf{M}^{N-1})_{x_1,x_0}$ which appears in (3.6) prevents that probability from simplifying even if $\nu$ is an invariant measure associated with $\mathbf{M}$.

We now give the main result about the convergence of cyclic Markov chains.

THEOREM 3. *Let $E$ be a finite state space and $\mathbf{V}$ be a square nonnegative matrix indexed by the elements of $E$ such that $\mathbf{V}$ admits a simple real eigenvalue $\lambda$ greatest in modulus than every other eigenvalues. Let $(X^{(N)})_{N \geq 1}$ be a family of stochastic processes such that for every $N$, $X^{(N)}$ is indexed by $\{0, \dots, N-1\}$ and*

$$(3.7) \qquad \mathbb{P}(X_0^{(N)} = x_0, \dots, X_{N-1}^{(N)} = x_{N-1}) = \frac{\prod_{i=0}^{N-1} \mathbf{V}_{x_i,x_{i+1}}}{\text{trace}(\mathbf{V}^N)},$$

*with the convention $x_N = x_0$.*

*Let $R = (R_i)_{i \in E}$ and $L = (L_i)_{i \in E}$ be respectively a right and a left eigenvector associated with $\lambda$ such that their dot product is equal to one, that is, $\sum L_i R_i = 1$.*

(i) *For each $N \geq 1$, $X^{(N)}$ is a $(\mathcal{U}_E, \mathbf{M}, N)$-cyclic MC, where $\mathbf{M}$ is equal to*

$$(3.8) \qquad \mathbf{M}_{i,j} = \mathbf{V}_{i,j} \frac{R_j}{\lambda \cdot R_i} \qquad \text{for } i, j \in E.$$

(ii) *Let now $X = (X_i)_{i \in \mathbb{N}}$ be a (well-defined) stochastic process and its finite-dimensional laws are given, for any $k \in \mathbb{N}$, by*

$$(3.9) \qquad \mu(\{x_0, \dots, x_k\}) = \lim_{N \to \infty} \mathbb{P}(X_0^{(N)} = x_0, \dots, X_k^{(N)} = x_k).$$

*Under $\mu$, $X$ is a $(\nu, \mathbf{M})$-MC, where $\mathbf{M}$ is defined as in equation (3.8) and $\nu$ is the invariant probability measure for $\mathbf{M}$ and is given by $\nu(x) = L_x R_x$, for $x \in E$.*

PROOF. We begin with (ii) and show that the limit in (3.9) exists. Let $k \in \mathbb{N}$ and $x_0, \ldots, x_k \in E$, for any $N > k$ we have

$$(3.10) \quad \mathbb{P}(X_0^{(N)} = x_0, \ldots, X_k^{(N)} = x_k) = \left( \prod_{i=0}^{k-1} \mathbf{V}_{x_i, x_{i+1}} \right) \frac{(\mathbf{V}^{N-k})_{x_k, x_0}}{\mathrm{trace}(\mathbf{V}^N)}.$$

When $N$ goes to infinity, the only significant terms of $(\mathbf{V}^{N-k})_{x_k, x_0}$ and $\mathrm{trace}(\mathbf{V}^N)$ are those in $\lambda^N$. More precisely,

$$(3.11) \qquad (\mathbf{V}^{N-k})_{x_k, x_0} = R_{x_k} L_{x_0} \lambda^{N-k} + \sum_{\lambda' \mathrm{eigenvalue\ of\ } \mathbf{V} \neq \lambda} a_{\lambda'} \lambda'^{N-k}$$

$$(3.12) \qquad = R_{x_k} L_{x_0} \lambda^{N-k} + o(\lambda^{N-k})$$

as $\lambda > |\lambda'|$, besides $\mathrm{trace}(V^N) = \lambda^N + o(\lambda^N)$ which leads to

$$(3.13) \qquad \lim_N \mathbb{P}(X_0^{(N)} = x_0, \ldots, X_k^{(N)} = x_k) = \frac{R_{x_k} L_{x_0}}{\lambda^k} \prod_{i=0}^{k-1} \mathbf{V}_{x_i, x_{i+1}}.$$

Let $\mu(\{x_0, \ldots, x_k\}) = \frac{R_{x_k} L_{x_0}}{\lambda^k} \prod_{i=0}^{k-1} \mathbf{V}_{x_i, x_{i+1}}$, we can check that $\nu$ is a probability distribution. Indeed, from equations (3.10) and (3.13)

$$\sum_{x \in E} \nu(x) = \sum_{x \in E} R_x L_x = \sum_{x \in E} \lim_N \frac{(\mathbf{V}^N)_{x,x}}{\mathrm{trace}(\mathbf{V}^N)} = \lim_N \sum_{x \in E} \frac{(\mathbf{V}^N)_{x,x}}{\mathrm{trace}(\mathbf{V}^N)} = 1,$$

where the inversion of the sum and the limit is immediate since $E$ is finite. We check similarly that the matrix $\mathbf{M}$ defined in (3.8) is stochastic.

Now it is easy to see that the finite-dimensional laws given in (3.9) are consistent, the Kolmogorov extension theorem applies and ensures that the stochastic process $X$ is well defined.

Point (i) follows directly from the definition of a cyclic Markov chain. $\square$

## 4. Examples of graphs.

We give in this section some examples of results that can be obtained by the application of Theorems 2 and 3. In the following examples, we only consider oriented lattices with vertices indexed by a subset of $\mathbb{Z}^2$. The $j$th line of the graph is the set of vertices with second coordinate equal to $j$. We will see why the restriction of the gas model "to a line" of the graph is Markovian. The general approach used is widely inspired by the method developed in [4] and [5]. For a given graph $G$, we first show that the assumptions of Theorem 2 are verified for $G$ and the sequence of lattices $(G^{(N)})_n$, that implies that the gas process on $G^{(N)}$ converges in distribution to the gas process on $G$. We then compute the distribution of the gas on a line of $G^{(N)}$ and interpret it as a cyclic MC by checking that its distribution can be written in a multiplicative form as in equation (3.7). Theorem 2 and Theorem 3 imply then that the gas process restricted to a line is Markovian. We explain fully the first example and sketch the others.

4.1. *The family of lattices* $(\mathcal{L}_R)_{R \subset \mathbb{N}}$. We define in this section a new family of lattices. For any finite subset $R$ of $\mathbb{N}$ such that $|R| \geq 2$, we define $\mathcal{L}_R$ as the lattice with set of vertices indexed by $\mathbb{Z}^2$ and from each vertex $(i, j)$, there are $|R|$ emerging edges from $(i, j)$ to $(i + r, j + 1)$ for $r \in R$. In the following, we always assume that $\inf(R) = 0$ without loss of generality. We set $\bar{R} = \sup(R)$.

Note that $\mathcal{L}_{\{0,1\}}$ corresponds to the square lattice. If $R = \{0, \ldots, n-1\}$, then $\mathcal{L}_R = \mathcal{L}_n$, which corresponds to the family of lattices introduced in [5] and detailed in the following. Another example is given in Figure 4.

REMARK 2. For any finite subset $R$ of $\mathbb{N}$, the lattice $\mathcal{L}_R$ verifies the assumption of Section 2 so the gas model is well defined for any $p < p_{\text{crit}}^{\mathcal{L}_R}$ and since the outdegree of any vertex is equal to $|R|$, $p_{\text{crit}}^{\mathcal{L}_R} > 1/|R| > 0$. For $N > n + \bar{R}$, the balls or radius $n$ of $\mathcal{L}_R^{(N)}$ and of $\mathcal{L}_R$ are isomorphic, moreover assumption 2 of Theorem 2 holds true with $d = |R|$, thus the finite-dimensional laws of the gas model on $\mathcal{L}_R$ converge to the ones on $\mathcal{L}_R^{(N)}$.

We denote $X_j^{(N)}$ the $N$-tuple that gives the occupation of the gas on the $j$th line of $\mathcal{L}_R^{(N)}$ and compute its distribution. The construction of the gas model given in Section 2 implies that, when $j$ decreases, $(X_j^{(N)})_{j \in \mathbb{Z}}$ is a "vertical" Markov chain (with $2^N$ states) under its stationary distribution. Markov chain theory implies that such a distribution is unique (that one of the main tool in [5] and [4]).

For $C \subset [N]$, let $F_C^{(N)}$ be the probability that the occupied vertices of a line of the graph are exactly those with first coordinate lying in $C$. In other

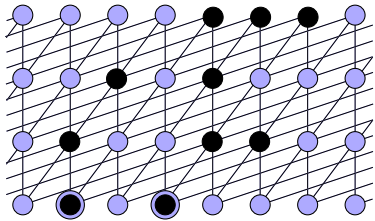

FIG. 4. *Example of a DA of size 10 on lattice $\mathcal{L}R$, when $R = \{0, 1, 4\}$.*

words, for a gas occupation $X^{\mathcal{L}_R}$ distributed according to the gas model given in Section 2:

(4.1) $$F_C^{(N)} = \mathbb{P}_p^{\mathcal{L}_R^{(N)}}(X^{\mathcal{L}_R^{(N)}}(i,j) = 1 \text{ if and only if } i \in C).$$

Note that the construction of the gas model implies that $F_C^{(N)}$ does not depend on a particular choice of $j$.

We traduce the fact that $(X_j^{(N)})_{j \in \mathbb{Z}, j \downarrow}$ is Markovian into recurrence relations for $F_C^{(N)}$. To that purpose, we define for any subset $C$ of $\mathbb{N}$

$$\mathcal{N}(C) = \bigcup_{i \in C} \{i + r \mid r \in R\}$$

and

$$\bar{\mathcal{N}}(C) = \bigcup_{i \in C} \{i - r \mid r \in R\},$$

where the addition is taken in $[N]$. Notice that $\{\mathcal{N}(C) \times \{1\}\}$ and $\{\bar{\mathcal{N}}(C) \times \{-1\}\}$ correspond respectively to the set of children and of fathers of the set $\{C \times \{0\}\}$ and that $|\mathcal{N}(C)| = |\bar{\mathcal{N}}(C)|$. We thus obtain the following equations:

(4.2) $$F_C^{(N)} = \left(\frac{p}{1-p}\right)^{|C|} \sum_{D \subset (\mathcal{N}(C))^c} (1-p)^{N - |\bar{\mathcal{N}}(D)|} F_D^{(N)}.$$

Following Bousquet-Mélou and Conway [5], we check that the probability distribution defined by

$$F_C^{(N)} = \frac{1}{Z_N} \left(\frac{p}{1-p}\right)^{|C|} (1-p)^{|\mathcal{N}(C)|},$$

where

$$Z_N = \sum_{C \subset [N]} p^{|C|} (1-p)^{\mathcal{N}(C) - |C|}$$

is stationary for equation (4.2).

To obtain a matrix formulation of that distribution, we consider the matrix $\mathbf{V}$ indexed by the elements of $\{0,1\}^{\bar{R}}$ and defined by, for $\sigma = (s_0, \ldots, s_{\bar{R}-1})$ and $\tau = (t_1, \ldots, t_{\bar{R}})$,

$$\mathbf{V}_{\sigma,\tau} = 0 \qquad \text{if } (s_1, \ldots, s_{\bar{R}-1}) \neq (t_1, \ldots, t_{\bar{R}-1});$$

and otherwise,

$$\mathbf{V}_{\sigma,\tau} = \begin{cases} p, & \text{if } t_{\bar{R}} = 1, \\ 1-p, & \text{if } t_{\bar{R}} = 0 \text{ and there exists } r \text{ such that } s_r = 1 \text{ and } \bar{R} - r \in R, \\ 1, & \text{otherwise.} \end{cases}$$

The quantity $F_C^{(N)}$ can then be rewritten as

$$F_C^{(N)} = \frac{1}{Z_N} \prod_{i=0}^{N-1} \mathbf{V}_{\sigma_i, \sigma_{i+1}},$$

(4.3)

$$\text{where } \begin{cases} \sigma_N = \sigma_0 & \text{and} \\ \sigma_i(k) = 1, & \text{iff } i+k-1 \in C, \end{cases}$$

with

$$Z_N = \sum_{\sigma_1, \ldots, \sigma_N} \left( \prod_{i=0}^{N-1} \mathbf{V}_{\sigma_i, \sigma_{i+1}} \right) = \text{trace}(\mathbf{V}^N).$$

The expression given for $F_C^{(N)}$ in equation (4.3) is in the very same form as the statement of Theorem 3. Furthermore it is immediate to check that all the coefficients of $\mathbf{V}^{\bar{R}-1}$ are positive [for $p \in (0,1]$] which ensures that $\mathbf{V}$ satisfies the conditions of Theorem 3 by the Perron–Frobenius theorem.

As mentioned in Remark 2, the finite dimensional laws of the gas occupation on $\mathcal{L}_R^{(N)}$ converges to those on $\mathcal{L}_R$. We apply Theorem 3 and get:

THEOREM 4. *Let $X = (X(i,j))_{(i,j) \in \mathbb{Z}^2}$ be the gas process on $\mathcal{L}_R$ distributed according to $\mathbb{P}_p^{\mathcal{L}_R}$, with $p < p_{\mathrm{crit}}^{\mathcal{L}_R}$. The stochastic process $(\boldsymbol{\Sigma}_i)_{i \in \mathbb{N}}$ defined by $\boldsymbol{\Sigma}_i = (X(i,0), \ldots, X(i+\bar{R}-1, 0))$ is a Markov chain under its stationary distribution.*

*In other words, the stochastic process $(X_i = X(i,0))_{i \in \mathbb{N}}$ is a Markov chain with memory $\bar{R}-1$ under its stationary distribution.*

*The family of lattices $\mathcal{L}_n$.* The family of lattices $(\mathcal{L}_n)_{n \geq 2}$ introduced in [5] corresponds to the particular case of $\mathcal{L}_R$ when $R = \{0, 1, \ldots, n-1\}$ (examples of DA on $\mathcal{L}_3$ and $\mathcal{L}_4$ are given in Figure 5). In [5] the GF for DA with one source is given as the solution of an algebraic solution of degree at most $n+1$.

THEOREM 5 (Bousquet-Mélou and Conway [5]). *The generating function $\mathcal{G}$ for DA on $\mathcal{L}_n$ with a single source is solution of the following equation:*

(4.4)    $t^2(1+t)^{n-1}[1+(n+1)\mathcal{G}]^{n+1} - [1+t+(n-1)\mathcal{G}]^{n-1}(t-2\mathcal{G}^2) = 0.$

We give some examples of computation obtained by the application of Theorem 4 on the lattices $\mathcal{L}_n$. In the case $n = 2$, the computation of the eigenvalues and eigenvectors of $\mathbf{V} = \begin{pmatrix} 1 & p \\ 1-p & p \end{pmatrix}$ constitutes an alternative proof of Theorem 3.3 of [10].

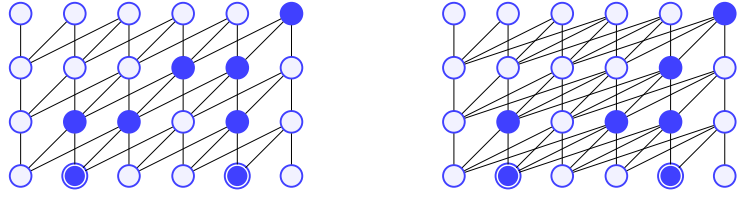

Fɪɢ. 5. *DA on the lattices $\mathcal{L}_3$ and $\mathcal{L}_4$ (all the edges are oriented upwards). The cells are black and the sources are circled.*

For $n = 3$, the transition matrix can be given explicitly as (the coefficients of the matrix are indexed by the lexicographical order on $\{0,1\}^2$)

$$
(4.5) \quad
\begin{pmatrix}
1/\lambda & 1-1/\lambda & 0 & 0 \\
0 & 0 & 1-p/2\lambda & p/2\lambda \\
\alpha & 1-\alpha & 0 & 0 \\
0 & 0 & 1-p/\lambda & p/\lambda
\end{pmatrix},
$$

$$
\text{where} \begin{cases}
\alpha = \dfrac{(1-p)^2 p}{(2-p-\lambda)\lambda}, \\[2mm]
\lambda = \dfrac{1+\sqrt{1+4p-4p^2}}{2}.
\end{cases}
$$

For example, we obtain as a consequence of that formula that the generating function $\mathcal{G}_k^{\mathcal{L}_3}$ for DA on $\mathcal{L}_3$ with a compact source of size $k \geq 2$ is equal to

$$
\mathcal{G}_k^{\mathcal{L}_3}(t) = \frac{1 - t(\sqrt{1-4t-4t^2})}{1-4t-4t^2+(1+2t)\sqrt{1-4t-4t^2}} \left( \frac{-2t}{1+\sqrt{1-4t-4t^2}} \right)^{k-1}.
$$

To obtain the GF as the solution of an algebraic equation, we use that in [5], the largest eigenvalue $\lambda$ of $\mathbf{V}$ is shown to be solution of

$$
(4.6) \quad \lambda^2(1-p)^{n-1} = \lambda^{n-1}(\lambda-1)^2.
$$

For $n \geq 4$, $\lambda$ cannot be computed explicitly from equation (4.6). Nevertheless, since $L$ and $R$ are eigenvectors associated to $\lambda$ their coordinates can be computed in linear time and are polynomial of degree one in $\lambda$. With the condition of renormalization $\sum L_i R_i = 1$, we obtain that for any free set $S$, the generating function for DA on $\mathcal{L}_n$ with source $S$ is a rational fraction and its numerator and denominator are polynomial in $\lambda$. Moreover, we know that $\lambda$ is solution of equation (4.6), which implies that the generating function is algebraic in $p$.

4.2. *The triangular lattice.* The triangular lattice, denoted Tri, is defined as the oriented graph with set of vertices $(i,j) \in \mathbb{Z}^2$ such that $i$ and $j$ have the same parity and with set of oriented edges $[(i,j),(i-1,j+1)]$, $[(i,j),(i+1,j+1)]$ and $[(i,j),(i,j+2)]$ (see Figure 6).

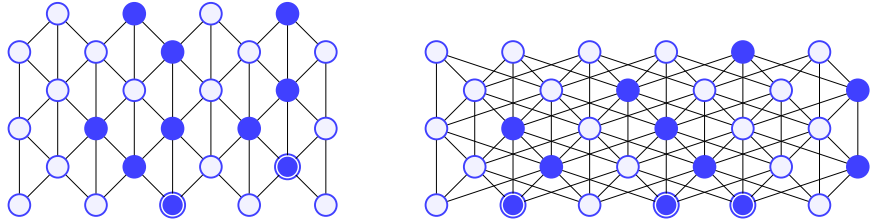

Fig. 6.   *Examples of DA on the triangular lattice (left) and on $\mathcal{T}_4$ (right).*

We follow some ideas used in [4] to compute the law of the gas on that lattice (note that in [4] the generating function for DA with one source on the triangular lattice is obtained by the study of an ad hoc gas model). We work on $\mathrm{Tri}^{(N)}$ the cylindric version of the triangular lattice. We keep the definition of $F_C^{(N)}$ introduced in equation (4.1), but since a vertex has children in the two following lines, we need to define an extension of $F_C^{(N)}$ to obtain recurrence relations. Let $C$ and $D$ be two subsets of $[N]$ and $X$ a gas model on $\mathrm{Tri}^{(N)}$, we denote $F_{C,D}^{(N)}$ the probability that the vertices occupied in the line 0 (respectively the line 1) of $\mathrm{Tri}^{(N)}$ are exactly the ones with first coordinate belonging to $C$ (respectively to $D$); in other words, for $C, D \in [N]$,

$$F_{C,D}^{(N)} = \mathbb{P}_p^{\mathrm{Tri}^{(N)}}\left( X^{\mathrm{Tri}^{(N)}}(i, \varepsilon) = 1 \text{ if and only if } \begin{cases} \varepsilon = 0 \text{ and } i/2 \in C \text{ or} \\ \varepsilon = 1 \text{ and } (i-1)/2 \in D \end{cases} \right).$$

We define $\mathcal{N}(C)$ as $\bigcup_{i \in C}\{i-1, i+1\}$ which leads to the following recurrence relation for $C, D \in [N]$ such that $\mathcal{N}(C) \cap D = \varnothing$:

$$(4.7) \qquad F_{C,D}^{(N)} = \left(\frac{p}{1-p}\right)^{|C|} \sum_{E \subset C^c(C \cup \mathcal{N}(D))} F_{D,E}(1-p)^{N-|\mathcal{N}(D) \cup E|}.$$

Notice that since the sum is taken on sets $E$ such that $\mathcal{N}(D) \cap E = \varnothing$, $|\mathcal{N}(D) \cup E|$ is equal to $|\mathcal{N}(D)| + |E|$. Therefore the distribution given by

$$(4.8) \qquad F_{C,D} = \frac{p^{|C|}p^{|D|}}{Z_N} 1_{\mathcal{N}(C) \cap D = \varnothing} \qquad \text{for } C, D \in [N],$$

and where

$$Z_N = \sum_{\substack{C,D \\ \mathcal{N}(C) \cap D = \varnothing}} p^{|C|}p^{|D|},$$

is solution to the recurrence relation given in (4.7).

Let $\mathbf{V} = \begin{pmatrix} 1 & p \\ 1 & 0 \end{pmatrix}$, we can rewrite equation (4.8) as

$$
(4.9) \quad F_{C,D} = \frac{1}{\text{trace}(\mathbf{V}^N)} \prod_{i=0}^{2N-1} \mathbf{V}_{x_i, x_{i+1}},
$$

$$
\text{where } x_i = \begin{cases} 1, & \text{if } i \text{ is even and } i \in C, \\ 1, & \text{if } i \text{ is odd and } i \in D, \\ 0, & \text{otherwise.} \end{cases}
$$

Combining equation (4.9) and Theorem 3 results in the following statement:

**Theorem 6.** *Let $X = (X(i,j))_{(i,j) \in \text{Tri}}$ be the gas process under $\mathbb{P}_p^{\text{Tri}}$, the stochastic process $\mathbf{\Sigma} = (\mathbf{\Sigma}_i)_{i \in \mathbb{Z}}$ defined by*

$$
\mathbf{\Sigma}_i = \begin{cases} X(i,0), & \text{if } i \text{ is even,} \\ X(i,1), & \text{if } i \text{ is odd,} \end{cases}
$$

*is a Markov chain under its stationary distribution and its transition matrix is given by*

$$
\mathbf{W} = \begin{pmatrix} \mathbb{P}(\mathbf{\Sigma}_1 = 0 | \mathbf{\Sigma}_0 = 0) & \mathbb{P}(\mathbf{\Sigma}_1 = 1 | \mathbf{\Sigma}_0 = 0) \\ \mathbb{P}(\mathbf{\Sigma}_1 = 0 | \mathbf{\Sigma}_0 = 1) & \mathbb{P}(\mathbf{\Sigma}_1 = 1 | \mathbf{\Sigma}_0 = 1) \end{pmatrix} = \begin{pmatrix} 1/\lambda & p/\lambda^2 \\ 1 & 0 \end{pmatrix},
$$

*where $\lambda = \frac{1 + \sqrt{1+4p}}{2}$ and its stationary distribution is given by*

$$
[\mathbb{P}(\mathbf{\Sigma}_0 = 0), \mathbb{P}(\mathbf{\Sigma}_0 = 1)] = [\lambda^2/(p+\lambda^2), p/(p+\lambda^2)].
$$

Adding up equation (4.8) for all possible $D$ leads to

$$
(4.10) \quad F_C = \frac{1}{Z_N} p^{|C|} (1+p)^{N - |\mathcal{N}(C)|}.
$$

Setting $\mathbf{V} = \begin{pmatrix} 1+p & p \\ 1 & p \end{pmatrix}$ enables equation (4.10) to be rewritten as

$$
(4.11) \quad F_C = \frac{1}{\text{trace}(V^N)} \prod_{i=0}^{N-1} V_{x_i, x_{i+1}}, \quad \text{where } x_i = 1 \text{ if and only if } 2i \in C.
$$

Again Theorem 3 and equation (4.11) lead to:

**Theorem 7.** *Let $X = (X(i,j))_{(i,j) \in \text{Tri}}$ be the gas process under $\mathbb{P}_p^{\text{Tri}}$, the stochastic process $\mathbf{\Sigma} = (\mathbf{\Sigma}_i)_{i \in \mathbb{Z}}$ defined by $\mathbf{\Sigma}_i = X(2i, 0)$ is a Markov chain under its stationary distribution and its transition matrix is given by*

$$
(4.12) \quad \begin{aligned} \mathbf{W} &= \begin{pmatrix} \mathbb{P}(\mathbf{\Sigma}_1 = 0 | \mathbf{\Sigma}_0 = 0) & \mathbb{P}(\mathbf{\Sigma}_1 = 1 | \mathbf{\Sigma}_0 = 0) \\ \mathbb{P}(\mathbf{\Sigma}_1 = 0 | \mathbf{\Sigma}_0 = 1) & \mathbb{P}(\mathbf{\Sigma}_1 = 1 | \mathbf{\Sigma}_0 = 1) \end{pmatrix} \\ &= \begin{pmatrix} 1 - \alpha_\circ & \alpha_\circ \\ \alpha_\bullet & 1 - \alpha_\bullet \end{pmatrix}, \end{aligned}
$$

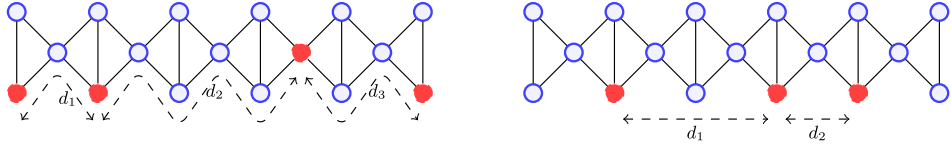

FIG. 7. *Examples of sources considered in Proposition 1(i) on the left and in Proposition 1(ii) on the right.*

and its stationary distribution by $[\mathbb{P}(\boldsymbol{\Sigma}_0 = 0), \mathbb{P}(\boldsymbol{\Sigma}_0 = 1)] = [\frac{\alpha_\bullet}{\alpha_\circ + \alpha_\bullet}, \frac{\alpha_\circ}{\alpha_\bullet + \alpha_\circ}]$, where

$$\alpha_\circ = \frac{2p}{1 + \sqrt{1 + 4p}} \quad and \quad \alpha_\bullet = \frac{1 + \sqrt{1 + 4p}}{2}.$$

The link between gas distribution and enumeration of DA given in Proposition 1 and a simple matrix computation give the following reinterpretation in terms of enumeration of DA of Theorems 6 and 7 (see Figure 7 for an example of the different sources considered).

PROPOSITION 1. (i) *Let $S = \{s_1, \ldots, s_k\}$ where $s_i = (x_i, \epsilon_i)$ be some points on the triangular lattice with $\epsilon_i \in \{0, 1\}$ and such that $d_i := x_{i+1} - x_i$ for $i \in \{1, \ldots, k-1\}$ are non smaller than 2. The GF of DA on the triangular lattice with source $S$ is given by*

$$\mathcal{G}_S^{\mathrm{Tri}}(-p) = (-1)^{|S|} \frac{\alpha}{1+\alpha} \prod_{i=1}^{k-1} \frac{(-\alpha)^{d_i} + \alpha}{1 + \alpha},$$

$$where\ \alpha = \frac{1 + 2p^2 + \sqrt{1 + 4p^2}}{2p^2}.$$

(ii) *Let $S = \{s_1, \ldots, s_k\}$ where $s_i = (2x_i, 0)$ be some vertices on a line of the triangular lattice, such that $d_i := x_{i+1} - x_i$ for $i \in \{1, \ldots, k-1\}$ are positive integers. The GF of DA on the triangular lattice with source $S$ is given by*

$$\mathcal{G}_S^{\mathrm{Tri}}(-p) = (-1)^{|S|} \frac{\alpha_\circ}{\alpha_\bullet + \alpha_\circ} \prod_{i=1}^{k-1} \frac{\alpha_\bullet(1 - \alpha_\bullet - \alpha_\circ)^{d_i} + \alpha_\circ}{\alpha_\bullet + \alpha_\circ}.$$

In particular, if $S_n := \{(i, 0), i = 1, \ldots, n\}$, we obtain $\mathcal{G}_{S_n}^{\mathrm{Tri}}(-p) = \frac{\alpha_\circ}{\alpha_\circ + \alpha_\bullet}(1 - \alpha_\bullet)^{n-1}(-1)^n$. Then, the GF of DA on the triangular lattice with compact sources satisfies

$$\sum_{n \geq 1} \mathcal{G}_{S_n}^{\mathrm{Tri}}(-p) = \sum_{n \geq 1} \frac{\alpha_\circ}{\alpha_\circ + \alpha_\bullet}(1 - \alpha_\bullet)^{n-1}(-1)^n = \frac{-p}{1 + 4p}.$$

That formula was obtained in [8] by combinatorial methods.

4.3. *The family of lattices $\mathcal{T}_n$.* We now study the family of lattices $\mathcal{T}_n$ introduced by Corteel, Denise and Gouyou-Beauchamps in [6]. The oriented lattice $\mathcal{T}_n$ is a combination of the lattice $\mathcal{L}_n$ and the triangular lattice, defined as follows:

- if $n = 2k + 1$, the vertices of $\mathcal{T}_n$ are labeled by the elements of $\mathbb{Z}^2$. From each vertex $(i, j) \in \mathbb{Z}^2$ there are $n$ emerging edges from $(i, j)$ to $(i + r, j + 1)$ for $-k \le r \le k$ and one emerging edge from $(i, j)$ to $(i, j + 2)$;
- if $n = 2k$, the vertices are labeled by the elements $(i, j) \in \mathbb{Z}^2$ such that $i$ and $j$ have the same parity. From each vertex $(i, j)$ there are $n$ emerging edges from $(i, j)$ to $(i + 2r + 1, y + 1)$ for $-k \le r \le k - 1$ and one emerging edge from $(i, j)$ to $(i, j + 2)$.

The case $n = 2$ corresponds to the triangular lattice, treated separately in Section 4.2 for sake of clarity. In [6], the generating function for DA on $\mathcal{T}_n$ with a single source is shown to be solution of an algebraic equation given explicitly. The proof relies on a combinatorial argument which links the generating function for DA on $\mathcal{T}_n$ to that for DA on $\mathcal{L}_n$.

The method used to obtain a stationary distribution for the gas model on $\mathcal{T}_n$ is very similar to that used in the case of the triangular lattice in Section 4.2. We keep the same definitions for $F_C^{(N)}$ and $F_{C,D}^{(N)}$ as those given for the triangular lattice and define for $C \in [N]$, $\mathcal{N}(C)$ as the set,

- $\bigcup_{i \in C} \{i + r, \text{ for } -k \le r \le k\}$ if $n = 2k + 1$,
- $\bigcup_{i \in C} \{i + 2r + 1, \text{ for } -k \le r \le k - 1\}$ if $n = 2k$.

With that new definition of $\mathcal{N}(C)$, equations (4.7), (4.8) and (4.10) still hold true. Following the ideas introduced to study $\mathcal{L}_n$ given in [5], we define $\mathbf{V}$ as the square matrix $(\mathbf{V}_{\sigma,\tau})_{\sigma,\tau}$ with indices running over $\{0, 1\}^{n-1}$ and defined as follows. If $\sigma = (s_1, \ldots, s_{n-1})$ and $\tau = (t_2, \ldots, t_n)$, then

$$(4.13) \quad \mathbf{V}_{\sigma,\tau} = \begin{cases} 0, & \text{if } (s_2, \ldots, s_{n-1}) \ne (t_2, \ldots, t_{n-1}), \\ p, & \text{if } (s_2, \ldots, s_{n-1}) = (t_2, \ldots, t_{n-1}) \text{ and } s_1 = 1, \\ 1 + p, & \text{if } \sigma = \tau = (0, 0, \ldots, 0), \\ 1, & \text{otherwise.} \end{cases}$$

The stationary distribution of the gas model on a line of $\mathcal{T}_n^{(N)}$ is given by

$$(4.14) \quad F_D = \frac{1}{Z_N} \prod_{i=0}^{N-1} \mathbf{V}_{\sigma_i, \sigma_{i+1}},$$

where $\begin{cases} \sigma_N = \sigma_0 & \text{and} \\ \sigma_i(k) = 1, & \text{if and only if } i + k - 1 \in D, \end{cases}$

with

$$Z_N = \sum_{\sigma_0, \ldots, \sigma_{N-1}} \left( \prod_{i=0}^{N-1} \mathcal{V}_{\sigma_i, \sigma_{i+1}} \right) = \text{trace}(\mathbf{V}^N).$$

The characteristic polynomial of $\mathbf{V}$, denoted $\chi$ can be calculated explicitly

$$\chi(x) = x^{2^{n-1}-n}\left(x^n - x^{n-1}(1+2p) + p^2\sum_{k=0}^{n-2} x^k\right).$$

We rewrite the latter equation as

$$\chi(x) = \frac{x^{2^{n-1}-n}}{1-x}(p^2 - x^{n-1}(x+p^2-1)(x-(2p+1))).$$

That implies that the dominant eigenvalue $\lambda$ of $\mathbf{V}$ satisfies $\lambda \neq 1$ and

$$(4.15) \qquad p^2 = \lambda^{n-1}(\lambda + p^2 - 1)(\lambda - (2p+1)).$$

We are here in the very same situation as for $\mathcal{L}_n$. We can compute explicitly the solutions of equation (4.15) only for $n < 4$. Nevertheless the same arguments as those given for $\mathcal{L}_n$ apply and we obtain the following from Theorem 3 and equation (4.14).

THEOREM 8.  *Let* $X = (X(i,j))_{(i,j)\in\mathcal{T}_n}$ *be the gas process under* $\mathbb{P}_p^{\mathcal{T}_n}$. *The stochastic process* $(\boldsymbol{\Sigma}_i)_{i\in\mathbb{N}}$ *defined by*

$$\boldsymbol{\Sigma}_i = (X(2i,0), X(2(i+1),0), \ldots, X(2(n-2),0)) \qquad \text{for } i \in \mathbb{N},$$

*is a Markov chain under its stationary distribution.*

*In other words,* $(X(2i,0))_{i\in\mathbb{N}}$ *is a Markov chain with memory* $n-1$ *under its stationary distribution.*

**Acknowledgment.**  I am very grateful to Jean-François Marckert who introduced me to this topic. This work has widely benefited from his various comments and suggestions.

LIAFA
UNIVERSITÉ PARIS DIDEROT—PARIS 7
CASE 7014
F-75205 PARIS CEDEX 13
FRANCE
E-MAIL: albenque@liafa.jussieu.fr
URL: http://www.liafa.jussieu.fr/~albenque